\documentclass[12pt,oneside,reqno]{amsart}
\usepackage{graphicx}

\graphicspath{ {figs/} }

\usepackage{pict2e}
\usepackage{amssymb}
\usepackage{amsthm}                
\usepackage[margin=1in]{geometry}  

\usepackage{epstopdf}
\usepackage{amsmath}
\usepackage{graphicx}
\usepackage{subfigure}
\usepackage{latexsym}
\usepackage{amsmath}
\usepackage{amsfonts}
\usepackage{amssymb}
\usepackage{tikz}
\usepackage{tikz-cd}
\usepackage{comment}
\usepackage{mathdots}
\usepackage{comment}
\usepackage{float}
\usepackage[mathscr]{euscript}
\usepackage{enumerate}
\usepackage{epsfig}
\usepackage { graphics }
\usepackage { graphicx }
\usepackage{listings}
\usepackage{algorithm}
\usepackage{fixltx2e}
\usepackage{algorithmic}
\usepackage{makecell}
\usepackage{nameref}
\usepackage{appendix}
\usepackage {eufrak}
\usepackage[utf8]{inputenc}
\usepackage{longtable}
\usepackage[lite]{amsrefs}

\renewcommand{\PrintDOI}[1]{\href{http://dx.doi.org/\detokenize{#1}}{doi: \detokenize{#1}}%
	\IfEmptyBibField{pages}{, (to appear in print)}{}}

\newtheorem{thm}{Theorem}[section]
\newtheorem{prop}[thm]{Proposition}

\theoremstyle{definition}
\newtheorem{theorem}{Theorem}[section]

\newtheorem{conjecture}[theorem]{Conjecture}

\theoremstyle{definition}
\newtheorem{definition}[theorem]{Definition}
\newtheorem{example}[theorem]{Example}
\theoremstyle{remark}

\numberwithin{equation}{section}

\numberwithin{equation}{section}

\usepackage[english]{babel}

\usepackage{amsmath}
\usepackage{graphicx}
\usepackage[colorlinks=true, allcolors=blue]{hyperref}
\title{State sum invariants of knots from idempotents in quandle rings}
\author[Elhamdadi]{Mohamed Elhamdadi} 
\address{Department of Mathematics and Statistics, University of South Florida, Tampa, FL, USA}
\email{emohamed@usf.edu}

\author[Swain]{Dipali Swain} 
\address{Department of Mathematics and Statistics, University of South Florida, Tampa, FL, USA}
\email{dipaliswain@usf.edu}

\begin{document}
\maketitle

\begin{abstract}

We use \emph{idempotents} in quandle rings in combination with the state sum invariants of knots to distinguish all of the $12965$ prime oriented knots up to $13$ crossings using \emph{only $21$ connected} quandles and three quandles made of idempotents in quandle rings.  We also distinguish all knots up to $13$ crossings from their mirror images using the same $24$ quandles. 
Furthermore, we distinguish all of the 2977 prime oriented knots up to 12 crossings using \emph{only} 10 connected quandles and three quandles made of idempotents in quandle rings.  This improves a result in [Quandle colorings of knots and applications, J. Knot Theory Ramifications 23 (2014), no. 6, 1450035].   
Our computations are achieved with the help of Python and Maple softwares.  
\end{abstract}


\section{Introduction}
Quandles are sets with self-distributive binary operations that axiomatize the three Reidemeister moves in classical knot theory \cites{EN, Joyce, Matveev}.  They have been used extensively to construct invariants of knots in the $3$-space and knotted surfaces in $4$-space, see for example \cites{CJKLS, CEGS}.  A cohomology theory for quandles was introduced in \cite{CJKLS} as a variation of rack homology defined by Fenn, Rourke and Sanderson in \cite{FRS}.  Two-dimensional cocycles were used to construct extensions of quandles \cite{CENS} and also allowed to define the quandle cocycle invariant of knots \cite{CJKLS}.  This homology was extended to a homology theory of set-theoretic Yang-Baxter equations in \cite{CES}.  

Quandle rings were introduced in \cite{BPS1} as a ``linearization" of quandles.  Their construction is similar to group rings \cite{Passman}.  Quandle rings have been investigated in \cite{EFT} where various properties of quandle rings were established such as decomposition of quandle rings into simple left or right ideals, the study of isomorphisms of quandle rings and a proof that  quandle rings are not power-associative rings.  In \cite{BPS2} the authors investigated the existence of nontrivial zero divisors in quandle rings in parallel to Kaplansky's study of zero divisors in group rings.  It was proved in \cite{ENSS} that integral quandle rings of quandles of finite type that are non-trivial coverings over nice base quandles admit infinitely many non-trivial idempotents, and gave their complete description.  In \cite{ENS} it was shown that quandle rings and their idempotents lead to proper enhancements of the well-known quandle coloring invariant of links. The authors also determined, for quandles of small cardinality, those for which the set of all idempotents is itself a quandle.  In \cite{CESY}, a set of $26$ finite quandles was used whose coloring invariant distinguishes all knots with $12$ or fewer crossings (modulo mirror images and reversal).  Furthermore, The authors of \cite{CESY} were able to find a set of $23$ quandles whose coloring invariant distinguishes $1058$ of these knots from their mirror images.  The results of both articles \cite{ENSS} and \cite{CESY} motivated our investigation.  Precisely, in this article, we combine the strength of \emph{idempotents} in quandle rings and the state sum invariants of knots to distinguish all $12965$ prime oriented knots up to $13$ crossings, using \emph{only $21$ connected} quandles and three quandles made of idempotents in quandle rings. We also distinguish all knots of up to $13$ crossings from their mirror images using the same \emph{only $21$ connected} quandles and three quandles made of idempotents in quandle rings.  As a byproduct, we distinguish all of the $2977$ prime oriented knots with up to $12$ crossings using \emph{only} $13$ quandles (10 connected quandles and three quandles made of idempotents in quandle rings), thus improving the result in \cite{CESY}.  Since the set of idempotents of a quandle ring may not be a quandle, 
we first use Maple software to select the quandles $X$ for which the set of idempotents $Y=\mathcal{I}(\mathbb{Z}_2[X])$ with $\mathbb{Z}_2$-coefficients is also a quandle (for classification of quandles up to isomorphism see \cite{EMR, VY, rig}).  We then consider a non-trivial $2$-cocycle $\phi$ of the quandle $X$ and a non-trivial $2$-cocycle $\psi$ of the quandle $Y=\mathcal{I}(\mathbb{Z}_2[X])$.  This allows the computations of the tuple $(\Phi_{(X, \phi)}(K), \Psi_{(Y, \psi)}(K))$ made of cocycle invariants of knots with respect to $(X,\phi)$ and $(Y,\psi)$.  Our computational results in comparison with the previous work in \cite{CESY} can be summarized in the following chart.\\

{\tiny{
\begin{tabular}{|c|c|c|c|}
     \hline
    Computation & Previous Work of \cite{CESY} & Our Work & Efficiency \\
     \hline
    Distinguishing knots up to 12 crossings & \# of quandles =  $26$  & \# of quandles= $13$ & $\approx 50\%$ \\
    \hline
    Distinguishing mirror image up to 12 crossings &\# of quandles =  $23$  & \# of quandles= $13$ & $\approx 56.52\%$ \\
    \hline
    Distinguishing knots up to $13$ crossings & - & \# of quandles =  $24$ & - \\
    \hline
     Distinguish mirror image up to $13$ crossings & - & \# of quandles= $24$ & - \\
     \hline
\end{tabular}
\newline
}\\
}

\par Furthermore, we distinguish all knots up to $13$ crossings from the list in \cite{CESY} using \emph{only $21$ connected} quandles and three quandles made of idempotents in quandle rings.

The article is organized as follows.  In Section~\ref{review} we review the basics of quandles, colorings and state sum invariants of knots.  Section~\ref{Idemp} gives the necessary material about idempotents of quandle rings.  In Section~\ref{Sec4}, for a quandle $X$ such that its idempotents over $\mathbb{Z}_2$ form a quandle $Y=\mathcal{I}(\mathbb{Z}_2[X])$, we consider a $2$-cocycle $\phi$ of $X$ with $\mathbb{Z}_2$-coefficients, a $2$-cocycle $\psi$ of $Y$ with $\mathbb{Z}_2$-coefficients and a 2-cocycle $\vartheta$ of a new quandle $W$ with $\mathbb{Z}_3$-coefficients.  We then compute the triples ($\Phi_{(X, \phi)}(K),\Psi_{(Y, \psi)}(K),\Theta_{(W, \vartheta)}(K)$) made of the cocycle invariants of a knot $K$ with respect to $(X,\phi)$, the cocycle invariants of a knot $K$ with respect to $(Y,\psi)$ and the cocycle invariants of a knot $K$ with respect to $(W,\vartheta)$.  The triples constructed here  make a strong invariant of knots.  In Section~\ref{Sec6} we distinguish knots with $13$ crossings in a similar way to the previous section except that here the quandle $W$ is obtained by \emph{iterating twice the idempotent process over $\mathbb{Z}_2$}, that is $W=\mathcal{I}(\mathbb{Z}_2[Y])= \mathcal{I}(\mathbb{Z}_2[\mathcal{I}(\mathbb{Z}_2[X])])$, here also we only consider the cases when both the set of idempotents $Y$ and $W$ are quandles.  Section~\ref{Sec7} gives a description of the algorithm.  

\section{Review of Quandles, Colorings and State Sum Invariants of Knots }\label{review}

In this section, we collect some basics about quandles, colorings and state sum invariants of knots needed for the rest of the article.  We begin with the following definition. 

\begin{definition}
    A quandle is a non empty set $X$ equipped with binary operation $\ast$ that satisfies the following axioms:
    \begin{enumerate}
        \item For every $x \in X,~ x\ast x=x$.
        \item For every $x,y \in X, \exists! z \in X$ such that $z\ast x=y$.
        \item For every $x,y,z \in X,~ (x\ast y)\ast z=(x\ast z)\ast(y\ast z)$.
    \end{enumerate}
   
\end{definition}
 \par
  A quandle homomorphism is a map between two quandles $f:(X,\ast) \to (Y,\star)$ such that $f(x\ast y)=f(x)\star f(y)$.  A quandle isomorphsim is a bijective quandle homomorphism. A quandle automorphsim is a quandle isomorphism of $X$ with itself. Each right multiplication $ S_y: X \to X$ given by $S_y(x)=x\ast y $ is an automorphism of $X$  fixing the element $y$. $S_y^{-1}(x)$ will be denoted by $x\overline{\ast }y=z \iff z\overline{\ast }y=x$.

\par A quandle is completely determined by the set of its right multiplications.  

\par Some typical examples of quandles are
\begin{itemize}
    \item Any set $X$ with operation $\ast$ given by $x\ast y=x$ for any $x,y \in X$ is a quandle known as \textit{trivial quandle}.
\item
A conjugacy class $X$ of a group $G$ with the quandle operation $x*y=y^{-1}xy $.
We call this a {\it conjugation quandle}.

\item
Any ${\mathbb{Z} }[T, T^{-1}]$-module $M$ is a quandle with
$x*y=Tx+(1-T)y$, $x,y \in M$, called an {\it  Alexander  quandle}.
\end{itemize}

The set ${\rm Aut}(X)$ of all automorphisms of a quandle $X$ is a group.  The subgroup ${\rm Inn}(X)$ generated by the set $\{S_x, \; x \in X\}$ is a normal subgroup of ${\rm Aut}(X)$ and is called the inner group of $X$.  When the action of ${\rm Inn}(X)$ on $X$ is transitive then the quandle $X$ is said to be a \emph{connected} quandle.  Through this article, we will consider connected quandles and some of their idempotents.  For more on quandles, the reader can consult \cite{EN, Joyce, Matveev}. 

\begin{figure} [ht]  
    
  \begin{center} 
\begin{tikzpicture}[x=0.75pt,y=0.75pt,yscale=-1,xscale=1]

\draw    (162.67,71.33) -- (192.4,101.47) ;
\draw    (202,112.67) -- (238.31,148.72) ;
\draw [shift={(239.73,150.13)}, rotate = 224.8] [color={rgb, 255:red, 0; green, 0; blue, 0 }  ][line width=0.75]    (10.93,-3.29) .. controls (6.95,-1.4) and (3.31,-0.3) .. (0,0) .. controls (3.31,0.3) and (6.95,1.4) .. (10.93,3.29)   ;
\draw    (232.67,73.67) -- (151.2,148.78) ;
\draw [shift={(149.73,150.13)}, rotate = 317.32] [color={rgb, 255:red, 0; green, 0; blue, 0 }  ][line width=0.75]    (10.93,-3.29) .. controls (6.95,-1.4) and (3.31,-0.3) .. (0,0) .. controls (3.31,0.3) and (6.95,1.4) .. (10.93,3.29)   ;
\draw    (424.4,102.8) -- (452.4,73.47) ;
\draw    (411.87,115.57) -- (378.31,150.92) ;
\draw [shift={(376.93,152.37)}, rotate = 313.51] [color={rgb, 255:red, 0; green, 0; blue, 0 }  ][line width=0.75]    (10.93,-3.29) .. controls (6.95,-1.4) and (3.31,-0.3) .. (0,0) .. controls (3.31,0.3) and (6.95,1.4) .. (10.93,3.29)   ;
\draw    (384,73.67) -- (453.72,149.99) ;
\draw [shift={(455.07,151.47)}, rotate = 227.59] [color={rgb, 255:red, 0; green, 0; blue, 0 }  ][line width=0.75]    (10.93,-3.29) .. controls (6.95,-1.4) and (3.31,-0.3) .. (0,0) .. controls (3.31,0.3) and (6.95,1.4) .. (10.93,3.29)   ;

\draw (145.33,55.07) node [anchor=north west][inner sep=0.75pt]  [font=\normalsize]  {$x$};
\draw (245.33,145.07) node [anchor=north west][inner sep=0.75pt]  [font=\normalsize]  {$x\ast y$};
\draw (131.33,146.07) node [anchor=north west][inner sep=0.75pt]    {$y$};
\draw (241.33,55.4) node [anchor=north west][inner sep=0.75pt]    {$y$};
\draw (462,145.4) node [anchor=north west][inner sep=0.75pt]    {$y$};
\draw (368.67,56.73) node [anchor=north west][inner sep=0.75pt]    {$y$};
\draw (459.33,55.73) node [anchor=north west][inner sep=0.75pt]  [font=\normalsize]  {$x$};
\draw (338,144.73) node [anchor=north west][inner sep=0.75pt]    {$x\overline{\ast } y$};

\end{tikzpicture}

    \caption{Rules of colorings at crossing}
    \label{fig:col}
    \end{center}
\end{figure}
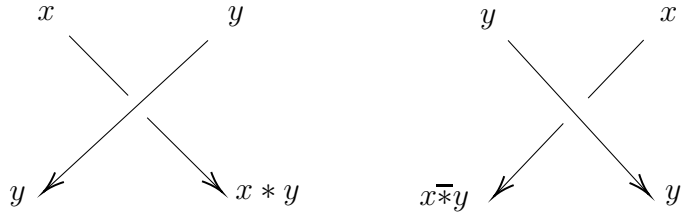

Let $X$ be a quandle and $K$ be a classical knot or a link diagram. Let $\mathcal{R}$ be a set of arcs. Then a coloring of $K$ by the quandle $X$ is a map from  $\mathcal{R}$ to $X$ such that at every crossing, the relations in Figure \ref{fig:col} hold. Alternately, it is a homomorphism from the knot quandle (see for example \cite{EN, Joyce, Matveev}) of $K$ to the quandle $X$.  It is well known that the number of colorings $\mathcal{C}_X(K)$ of a knot $K$ by a quandle $X$ is a knot invariant (see for example \cite{MR3381331}). 
 
Let $X$ be a finite quandle and let $A$ be an abelian group. Let $C_n^{\rm R}(X)$ be the free abelian group generated by $n$-tuples $(x_1, \dots, x_n)$ of elements $X$. We define a homomorphism $\partial_{n}: C_{n}^{\rm R}(X) \to C_{n-1}^{\rm R}(X)$ by 
\begin{eqnarray*}
\lefteqn{
\partial_{n}(x_1, x_2, \dots, x_n) } \nonumber \\ && =
\sum_{i=2}^{n} (-1)^{i}\left[ (x_1, x_2, \dots, x_{i-1}, x_{i+1},\dots, x_n) \right.
\nonumber \\
&&
- \left. (x_1 \ast x_i, x_2 \ast x_i, \dots, x_{i-1}\ast x_i, x_{i+1}, \dots, x_n) \right]
\end{eqnarray*}
for $n \geq 2$ 
and $\partial_n=0$ for 
$n \leq 1$. 
 Then
$C_\ast^{\rm R}(X)
= \{C_n^{\rm R}(X), \partial_n \}$ is a chain complex.

Let $C_n^{\rm D}(X)$ be the subset of $C_n^{\rm R}(X)$ generated
by $n$-tuples $(x_1, \dots, x_n)$
with $x_{i}=x_{i+1}$ for some $i \in \{1, \dots,n-1\}$ if $n \geq 2$;
otherwise let $C_n^{\rm D}(X)=0$. If $X$ is a quandle, then
$\partial_n(C_n^{\rm D}(X)) \subset C_{n-1}^{\rm D}(X)$ and
$C_\ast^{\rm D}(X) = \{ C_n^{\rm D}(X), \partial_n \}$ is a sub-complex of
$C_\ast^{\rm
R}(X)$. Consider the quotient complex $\{C_\ast^{\rm Q}(X)\}$ with   $C_n^{\rm Q}(X) = C_n^{\rm R}(X)/ C_n^{\rm D}(X)$.  
For quandles, the chain and cochain complexes with coefficient in an abelian group $A$ are given by
\begin{eqnarray*}
C_\ast^{\rm Q}(X;A) = C_\ast^{\rm Q}(X) \otimes A, \quad && \partial =
\partial \otimes {\rm id}; \\ C^\ast_{\rm Q}(X;A) = {\rm Hom}(C_\ast^{\rm
Q}(X), A), \quad
&& \delta= {\rm Hom}(\partial, {\rm id}).
\end{eqnarray*}

The $n$\/th {\it quandle homology group\/}  and the $n$\/th
{\it quandle cohomology group\/ } \cite{CJKLS} of a quandle $X$ with coefficient group $A$ are given by
\begin{eqnarray*}
H_n^{\rm Q}(X; A) 
 = H_{n}(C_\ast^{\rm Q}(X;A)), \quad
H^n_{\rm Q}(X; A) 
 = H^{n}(C^\ast_{\rm Q}(X;A)). \end{eqnarray*}
For more details on quandle cohomology see \cite{CJKLS}.  In this article we will focus on low dimensional cohomology and precisely $2$-cocycles as they are needed to define the quandle cocycle invariant of knots.
A function $\phi: X \times X \to A$ is called a quandle 2-cocycle if it satisfies the 2-cocycle condition: 
\begin{equation}
\label{cocy1}
    \phi(x,y)-\phi(x,z)+\phi(x\ast y, z)-\phi(x\ast z, y\ast z)=0;~~~~ \forall x, y, z \in X 
    \end{equation}
 and   
    \begin{equation}
    \label{cocy2}
      \phi(x,x)=0,~~~ \forall x \in X  
    \end{equation}
    
    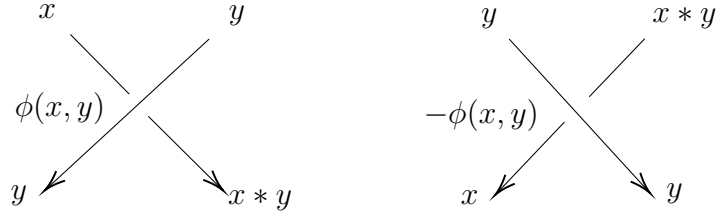
\begin{figure}[ht]
    \centering
\begin{tikzpicture}[x=0.75pt,y=0.75pt,yscale=-1,xscale=1]

\draw    (162.67,71.33) -- (192.4,101.47) ;
\draw    (202,112.67) -- (238.31,148.72) ;
\draw [shift={(239.73,150.13)}, rotate = 224.8] [color={rgb, 255:red, 0; green, 0; blue, 0 }  ][line width=0.75]    (10.93,-3.29) .. controls (6.95,-1.4) and (3.31,-0.3) .. (0,0) .. controls (3.31,0.3) and (6.95,1.4) .. (10.93,3.29)   ;
\draw    (232.67,73.67) -- (151.2,148.78) ;
\draw [shift={(149.73,150.13)}, rotate = 317.32] [color={rgb, 255:red, 0; green, 0; blue, 0 }  ][line width=0.75]    (10.93,-3.29) .. controls (6.95,-1.4) and (3.31,-0.3) .. (0,0) .. controls (3.31,0.3) and (6.95,1.4) .. (10.93,3.29)   ;
\draw    (424.4,102.8) -- (452.4,73.47) ;
\draw    (411.87,115.57) -- (378.31,150.92) ;
\draw [shift={(376.93,152.37)}, rotate = 313.51] [color={rgb, 255:red, 0; green, 0; blue, 0 }  ][line width=0.75]    (10.93,-3.29) .. controls (6.95,-1.4) and (3.31,-0.3) .. (0,0) .. controls (3.31,0.3) and (6.95,1.4) .. (10.93,3.29)   ;
\draw    (384,73.67) -- (453.72,149.99) ;
\draw [shift={(455.07,151.47)}, rotate = 227.59] [color={rgb, 255:red, 0; green, 0; blue, 0 }  ][line width=0.75]    (10.93,-3.29) .. controls (6.95,-1.4) and (3.31,-0.3) .. (0,0) .. controls (3.31,0.3) and (6.95,1.4) .. (10.93,3.29)   ;

\draw (145.33,55.07) node [anchor=north west][inner sep=0.75pt]  [font=\normalsize]  {$x$};
\draw (133.33,99.07) node [anchor=north west][inner sep=0.75pt]  [font=\normalsize]  {$\phi ( x,y)$};
\draw (131.33,146.07) node [anchor=north west][inner sep=0.75pt]    {$y$};
\draw (241.33,55.4) node [anchor=north west][inner sep=0.75pt]    {$y$};
\draw (462,145.4) node [anchor=north west][inner sep=0.75pt]    {$y$};
\draw (368.67,56.73) node [anchor=north west][inner sep=0.75pt]    {$y$};
\draw (358.33,148.73) node [anchor=north west][inner sep=0.75pt]  [font=\normalsize]  {$x$};
\draw (339.67,102.73) node [anchor=north west][inner sep=0.75pt]  [font=\normalsize]  {$-\phi ( x,y)$};
\draw (241,148.4) node [anchor=north west][inner sep=0.75pt]    {$x\ast y$};
\draw (455,56.4) node [anchor=north west][inner sep=0.75pt]    {$x\ast y$};

\end{tikzpicture}

    \caption{Boltzmann weights at crossing}
    \label{boltz}

\end{figure}


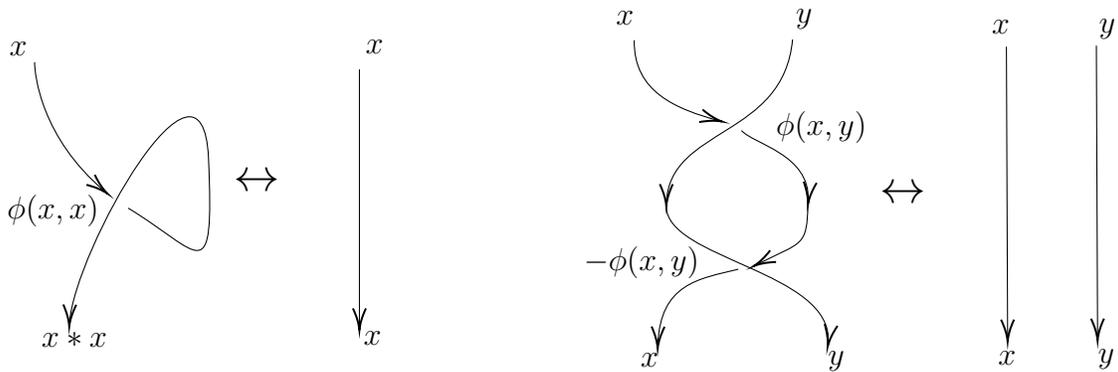
\begin{figure}[ht]
    \centering
\begin{tikzpicture}[x=0.75pt,y=0.75pt,yscale=-1,xscale=1]

\draw    (455.17,40.98) .. controls (450.68,90.87) and (393.26,84.97) .. (391.53,123.19) ;
\draw [shift={(391.49,124.96)}, rotate = 270] [color={rgb, 255:red, 0; green, 0; blue, 0 }  ][line width=0.75]    (10.93,-3.29) .. controls (6.95,-1.4) and (3.31,-0.3) .. (0,0) .. controls (3.31,0.3) and (6.95,1.4) .. (10.93,3.29)   ;
\draw    (375.34,41.38) .. controls (374.73,67.92) and (397.52,76.34) .. (417.61,81.91) ;
\draw [shift={(419.46,82.41)}, rotate = 195.2] [color={rgb, 255:red, 0; green, 0; blue, 0 }  ][line width=0.75]    (10.93,-3.29) .. controls (6.95,-1.4) and (3.31,-0.3) .. (0,0) .. controls (3.31,0.3) and (6.95,1.4) .. (10.93,3.29)   ;
\draw    (429.4,86.97) .. controls (435.46,93.89) and (463.97,98.88) .. (463.09,122.13) ;
\draw [shift={(462.96,123.95)}, rotate = 275.72] [color={rgb, 255:red, 0; green, 0; blue, 0 }  ][line width=0.75]    (10.93,-3.29) .. controls (6.95,-1.4) and (3.31,-0.3) .. (0,0) .. controls (3.31,0.3) and (6.95,1.4) .. (10.93,3.29)   ;
\draw    (391.49,124.96) .. controls (390.88,151.91) and (476.16,157.89) .. (473.11,194.2) ;
\draw [shift={(472.91,195.88)}, rotate = 278.82] [color={rgb, 255:red, 0; green, 0; blue, 0 }  ][line width=0.75]    (10.93,-3.29) .. controls (6.95,-1.4) and (3.31,-0.3) .. (0,0) .. controls (3.31,0.3) and (6.95,1.4) .. (10.93,3.29)   ;
\draw    (462.96,123.95) .. controls (462.35,145.79) and (460.55,142.83) .. (437.08,154.62) ;
\draw [shift={(435.62,155.36)}, rotate = 333] [color={rgb, 255:red, 0; green, 0; blue, 0 }  ][line width=0.75]    (10.93,-3.29) .. controls (6.95,-1.4) and (3.31,-0.3) .. (0,0) .. controls (3.31,0.3) and (6.95,1.4) .. (10.93,3.29)   ;
\draw    (427.83,156.98) .. controls (407.45,161.84) and (388.54,164.37) .. (387.2,194.98) ;
\draw [shift={(387.14,196.89)}, rotate = 271.1] [color={rgb, 255:red, 0; green, 0; blue, 0 }  ][line width=0.75]    (10.93,-3.29) .. controls (6.95,-1.4) and (3.31,-0.3) .. (0,0) .. controls (3.31,0.3) and (6.95,1.4) .. (10.93,3.29)   ;
\draw    (563.02,44.63) -- (563.63,191.35) ;
\draw [shift={(563.64,193.35)}, rotate = 269.76] [color={rgb, 255:red, 0; green, 0; blue, 0 }  ][line width=0.75]    (10.93,-3.29) .. controls (6.95,-1.4) and (3.31,-0.3) .. (0,0) .. controls (3.31,0.3) and (6.95,1.4) .. (10.93,3.29)   ;
\draw    (608.39,44.12) -- (609,190.84) ;
\draw [shift={(609.01,192.84)}, rotate = 269.76] [color={rgb, 255:red, 0; green, 0; blue, 0 }  ][line width=0.75]    (10.93,-3.29) .. controls (6.95,-1.4) and (3.31,-0.3) .. (0,0) .. controls (3.31,0.3) and (6.95,1.4) .. (10.93,3.29)   ;
\draw    (120.08,126.13) .. controls (156.78,148.23) and (163.57,169.38) .. (160.61,100.99) .. controls (157.7,33.62) and (93.23,146.94) .. (90.46,184.14) ;
\draw [shift={(90.37,185.79)}, rotate = 271.42] [color={rgb, 255:red, 0; green, 0; blue, 0 }  ][line width=0.75]    (10.93,-3.29) .. controls (6.95,-1.4) and (3.31,-0.3) .. (0,0) .. controls (3.31,0.3) and (6.95,1.4) .. (10.93,3.29)   ;
\draw    (72.9,52.48) .. controls (72.9,59.01) and (76.29,90.55) .. (108.12,117.54) ;
\draw [shift={(109.6,118.77)}, rotate = 219.38] [color={rgb, 255:red, 0; green, 0; blue, 0 }  ][line width=0.75]    (10.93,-3.29) .. controls (6.95,-1.4) and (3.31,-0.3) .. (0,0) .. controls (3.31,0.3) and (6.95,1.4) .. (10.93,3.29)   ;
\draw    (236.72,56.13) -- (236.72,185.67) ;
\draw [shift={(236.72,187.67)}, rotate = 270] [color={rgb, 255:red, 0; green, 0; blue, 0 }  ][line width=0.75]    (10.93,-3.29) .. controls (6.95,-1.4) and (3.31,-0.3) .. (0,0) .. controls (3.31,0.3) and (6.95,1.4) .. (10.93,3.29)   ;

\draw (445.59,76.03) node [anchor=north west][inner sep=0.75pt]    {$\phi ( x,y)$};
\draw (348.68,143.95) node [anchor=north west][inner sep=0.75pt]    {$-\phi ( x,y)$};
\draw (364.99,25.44) node [anchor=north west][inner sep=0.75pt]    {$x$};
\draw (455.72,24.43) node [anchor=north west][inner sep=0.75pt]    {$y$};
\draw (471.88,196.85) node [anchor=north west][inner sep=0.75pt]    {$y$};
\draw (377.42,198.07) node [anchor=north west][inner sep=0.75pt]    {$x$};
\draw (557.64,197.03) node [anchor=north west][inner sep=0.75pt]    {$x$};
\draw (554.43,29.96) node [anchor=north west][inner sep=0.75pt]    {$x$};
\draw (607.98,195.81) node [anchor=north west][inner sep=0.75pt]    {$y$};
\draw (608.6,29.26) node [anchor=north west][inner sep=0.75pt]    {$y$};
\draw (57.87,118.09) node [anchor=north west][inner sep=0.75pt]    {$\phi ( x,x)$};
\draw (59.04,40.76) node [anchor=north west][inner sep=0.75pt]    {$x$};
\draw (75.25,187.32) node [anchor=north west][inner sep=0.75pt]    {$x\ast x$};
\draw (238.59,39.69) node [anchor=north west][inner sep=0.75pt]    {$x$};
\draw (237.71,186.7) node [anchor=north west][inner sep=0.75pt]    {$x$};
\draw (172.01,107.4) node [anchor=north west][inner sep=0.75pt]  [font=\Large]  {$\leftrightarrow $};
\draw (498.01,113.4) node [anchor=north west][inner sep=0.75pt]  [font=\Large]  {$\leftrightarrow $};

\end{tikzpicture}

    \caption{Boltzmann weights from Reidmeister I and II}
    \label{fig:cocycle1}
\end{figure}


  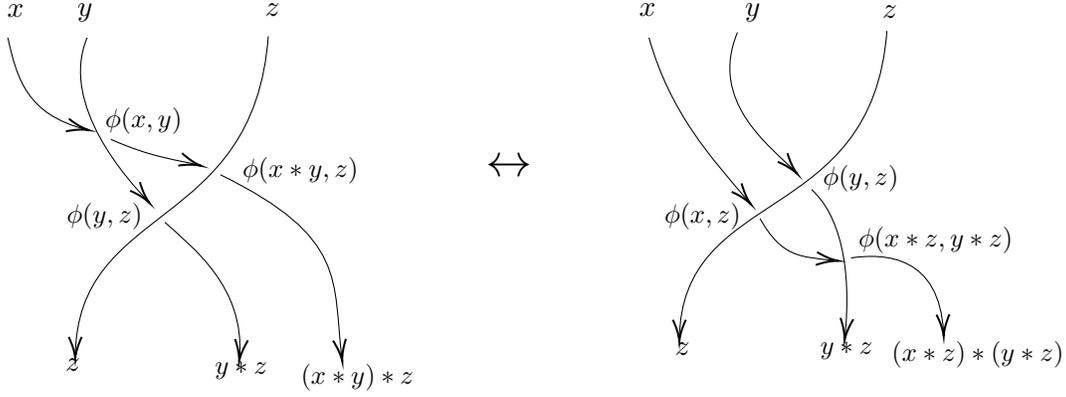
\begin{figure}[ht]
      \centering
\begin{tikzpicture}[x=0.75pt,y=0.75pt,yscale=-1,xscale=1]

\draw    (96.13,25.47) .. controls (81.72,59.93) and (117.97,98.52) .. (126.24,107.46) ;
\draw [shift={(127.47,108.8)}, rotate = 228.81] [color={rgb, 255:red, 0; green, 0; blue, 0 }  ][line width=0.75]    (10.93,-3.29) .. controls (6.95,-1.4) and (3.31,-0.3) .. (0,0) .. controls (3.31,0.3) and (6.95,1.4) .. (10.93,3.29)   ;
\draw    (187.47,24.8) .. controls (180.17,124.3) and (91.03,108.3) .. (90.14,186.28) ;
\draw [shift={(90.13,187.47)}, rotate = 270] [color={rgb, 255:red, 0; green, 0; blue, 0 }  ][line width=0.75]    (10.93,-3.29) .. controls (6.95,-1.4) and (3.31,-0.3) .. (0,0) .. controls (3.31,0.3) and (6.95,1.4) .. (10.93,3.29)   ;
\draw    (135.33,118.67) .. controls (149.25,131.27) and (175.98,151.52) .. (172.96,187.16) ;
\draw [shift={(172.8,188.8)}, rotate = 276.23] [color={rgb, 255:red, 0; green, 0; blue, 0 }  ][line width=0.75]    (10.93,-3.29) .. controls (6.95,-1.4) and (3.31,-0.3) .. (0,0) .. controls (3.31,0.3) and (6.95,1.4) .. (10.93,3.29)   ;
\draw    (56.13,25.47) .. controls (60.71,46.37) and (65.92,64.08) .. (94.99,71.68) ;
\draw [shift={(96.8,72.13)}, rotate = 193.45] [color={rgb, 255:red, 0; green, 0; blue, 0 }  ][line width=0.75]    (10.93,-3.29) .. controls (6.95,-1.4) and (3.31,-0.3) .. (0,0) .. controls (3.31,0.3) and (6.95,1.4) .. (10.93,3.29)   ;
\draw    (163.33,94.67) .. controls (227.82,126.32) and (219.86,148.27) .. (224.58,187.02) ;
\draw [shift={(224.8,188.8)}, rotate = 262.57] [color={rgb, 255:red, 0; green, 0; blue, 0 }  ][line width=0.75]    (10.93,-3.29) .. controls (6.95,-1.4) and (3.31,-0.3) .. (0,0) .. controls (3.31,0.3) and (6.95,1.4) .. (10.93,3.29)   ;
\draw    (499.47,22.13) .. controls (492.17,121.63) and (395.77,99.69) .. (394.81,177.62) ;
\draw [shift={(394.8,178.8)}, rotate = 270] [color={rgb, 255:red, 0; green, 0; blue, 0 }  ][line width=0.75]    (10.93,-3.29) .. controls (6.95,-1.4) and (3.31,-0.3) .. (0,0) .. controls (3.31,0.3) and (6.95,1.4) .. (10.93,3.29)   ;
\draw    (424.13,22.8) .. controls (409.34,58.18) and (441.1,80.53) .. (452.93,92.83) ;
\draw [shift={(454.13,94.13)}, rotate = 228.37] [color={rgb, 255:red, 0; green, 0; blue, 0 }  ][line width=0.75]    (10.93,-3.29) .. controls (6.95,-1.4) and (3.31,-0.3) .. (0,0) .. controls (3.31,0.3) and (6.95,1.4) .. (10.93,3.29)   ;
\draw    (461.47,102.13) .. controls (475.39,114.74) and (481.94,140.67) .. (478.31,176.49) ;
\draw [shift={(478.13,178.13)}, rotate = 276.23] [color={rgb, 255:red, 0; green, 0; blue, 0 }  ][line width=0.75]    (10.93,-3.29) .. controls (6.95,-1.4) and (3.31,-0.3) .. (0,0) .. controls (3.31,0.3) and (6.95,1.4) .. (10.93,3.29)   ;
\draw    (379.47,25.47) .. controls (390.05,63.16) and (410.62,85.88) .. (429.64,108.11) ;
\draw [shift={(430.8,109.47)}, rotate = 229.54] [color={rgb, 255:red, 0; green, 0; blue, 0 }  ][line width=0.75]    (10.93,-3.29) .. controls (6.95,-1.4) and (3.31,-0.3) .. (0,0) .. controls (3.31,0.3) and (6.95,1.4) .. (10.93,3.29)   ;
\draw    (435.47,116.8) .. controls (446.4,136.1) and (455.48,136.15) .. (472.87,137.93) ;
\draw [shift={(474.8,138.13)}, rotate = 186.12] [color={rgb, 255:red, 0; green, 0; blue, 0 }  ][line width=0.75]    (10.93,-3.29) .. controls (6.95,-1.4) and (3.31,-0.3) .. (0,0) .. controls (3.31,0.3) and (6.95,1.4) .. (10.93,3.29)   ;
\draw    (481.33,136.67) .. controls (498.54,134.17) and (525.96,135.43) .. (528.06,174.33) ;
\draw [shift={(528.13,176.13)}, rotate = 268.12] [color={rgb, 255:red, 0; green, 0; blue, 0 }  ][line width=0.75]    (10.93,-3.29) .. controls (6.95,-1.4) and (3.31,-0.3) .. (0,0) .. controls (3.31,0.3) and (6.95,1.4) .. (10.93,3.29)   ;
\draw    (108.13,76.8) .. controls (128.59,86.1) and (142.71,87.33) .. (151.56,89.6) ;
\draw [shift={(153.47,90.13)}, rotate = 197.1] [color={rgb, 255:red, 0; green, 0; blue, 0 }  ][line width=0.75]    (10.93,-3.29) .. controls (6.95,-1.4) and (3.31,-0.3) .. (0,0) .. controls (3.31,0.3) and (6.95,1.4) .. (10.93,3.29)   ;

\draw (54.67,7.4) node [anchor=north west][inner sep=0.75pt]  [font=\small]  {$x$};
\draw (184.67,6.73) node [anchor=north west][inner sep=0.75pt]  [font=\small]  {$z$};
\draw (84,186.07) node [anchor=north west][inner sep=0.75pt]  [font=\small]  {$z$};
\draw (90,6.73) node [anchor=north west][inner sep=0.75pt]  [font=\small]  {$y$};
\draw (84.67,107.4) node [anchor=north west][inner sep=0.75pt]  [font=\footnotesize]  {$\phi ( y,z)$};
\draw (173.33,84.4) node [anchor=north west][inner sep=0.75pt]  [font=\footnotesize]  {$\phi ( x\ast y,z)$};
\draw (391.33,178.07) node [anchor=north west][inner sep=0.75pt]  [font=\small]  {$z$};
\draw (496,8.07) node [anchor=north west][inner sep=0.75pt]  [font=\small]  {$z$};
\draw (426.67,6.73) node [anchor=north west][inner sep=0.75pt]  [font=\small]  {$y$};
\draw (373.33,6.73) node [anchor=north west][inner sep=0.75pt]  [font=\small]  {$x$};
\draw (484,119.73) node [anchor=north west][inner sep=0.75pt]  [font=\footnotesize]  {$\phi ( x\ast z,y\ast z)$};
\draw (466,88.73) node [anchor=north west][inner sep=0.75pt]  [font=\footnotesize]  {$\phi ( y,z)$};
\draw (104,59.4) node [anchor=north west][inner sep=0.75pt]  [font=\footnotesize]  {$\phi ( x,y)$};
\draw (159.33,187.73) node [anchor=north west][inner sep=0.75pt]  [font=\footnotesize]  {$y\ast z$};
\draw (202.67,188.73) node [anchor=north west][inner sep=0.75pt]  [font=\footnotesize]  {$( x\ast y) \ast z$};
\draw (464.67,177.73) node [anchor=north west][inner sep=0.75pt]  [font=\footnotesize]  {$y\ast z$};
\draw (500.67,177.4) node [anchor=north west][inner sep=0.75pt]  [font=\footnotesize]  {$( x\ast z) \ast ( y\ast z)$};
\draw (386,107.4) node [anchor=north west][inner sep=0.75pt]  [font=\footnotesize]  {$\phi ( x,z)$};
\draw (296,85.4) node [anchor=north west][inner sep=0.75pt]  [font=\Large]  {$\leftrightarrow $};

\end{tikzpicture}

       \caption{The quandle 2-cocycle condition (\ref{cocy1})  from Reidmeister III}
        \label{fig:cocycle3}
  \end{figure}

\par Let $X$ be a quandle and $\phi: X \times X \to A$ be a 2-cocycle. Consider a knot $K$ and let $\mathcal{C}_X(K)$ be a coloring of $K$. The Boltzmann weight at a crossing $\tau$ is defined by $\phi(x,y)^{\epsilon}$, where $\epsilon$ is the sign of the crossings (see Figure~\ref{boltz}).  Thus one sees that equation~(\ref{cocy1}) can be obtained from Figure~\ref{fig:cocycle3}.
\begin{definition}\cites{CJKLS}
    Let X be a quandle and $\phi$ be a 2-cocycle with coefficient in an abelian group $A$.  Let $D$ be a diagram of a knot $K$.   The state sum of the knot diagram $D$ is given by
    \begin{equation}
       \Phi(D)= \sum_{\mathcal{C}} \prod_\tau \phi(x,y)^{\epsilon}
    \end{equation}
    where the product is taken over all crossings of $D$ and the sum is taken over all the possible colorings of $D$.
\end{definition}

In \cites{CJKLS} the authors proved that the state sum is invariant under the three Reidemeister moves.

\begin{theorem}\cite[Theorem 4.4, page 3954]{CJKLS} 
    Let $\phi$ be a $2$-cocycle with coefficient in an abelian group $A$.  Let $D$ be a diagram of a knot $K$.  The state sum $\Phi(D)$ is invariant under the three Reidemeister moves, thus it is denoted by $\Phi_{\phi}(K)$.  
\end{theorem}
\section{Review of idempotents in quandle rings}\label{Idemp}

Given a \emph{non-trivial} quandle $(X, \ast)$ and an associative ring $\mathbf{k}$ with unity, we consider a \textit{nonassociative} ring $\mathbf{k}[X]$ called the quandle ring.   Precisely, Let $\mathbf{k}[X]$ be the set of elements that are uniquely expressible in the form $\sum_{x \in X }  a_x e_x$, where $x \in X$ and $a_x=0$ for almost all $x$.  Then the  set $\mathbf{k}[X]$ becomes a ring with the natural addition and the multiplication given by the following rule, where $x, y \in X$ and $a_x, b_y \in \mathbf{k}$,
\begin{eqnarray} \label{mult}
 ( \sum_{x \in X }  a_x e_x) \cdot ( \sum_{ y \in X }  b_y e_y )
=   \sum_{x, y \in X } a_x b_y e_{x \ast y} . \end{eqnarray}
For more on quandle rings, see for example \cite{BPS1, EFT}.
\begin{definition}
     Let $X$ be a quandle and $\mathbf{k}$ an integral domain with unity.  A non-zero element $u \in \mathbf{k}[X]$ is called an {\it idempotent} if $u^2=u$.  We denote the set of all idempotents of $\mathbf{k}[X]$ by $\mathcal{I}\big(\mathbf{k}[X]\big)$, that is,
$$
\mathcal{I}\big(\mathbf{k}[X] \big)=\big\{ v  \in \mathbf{k}[X],\;  \; v^2=v\big \}.
$$
\end{definition}
 From the equation~(\ref{mult}) it follows that the basis elements $\{e_x, x \in X\}$ are idempotents in  $\mathbf{k}[X]$.  We call them the {\it trivial idempotents}.
A non-trivial idempotent  is an element of $\mathbf{k}[X]$ that is not of the form $e_x$ for any $x \in X$.  Furthermore, for any quandle $X$ and  $Y$ such that $X \subset Y$, we have $\mathcal{I}\big(\mathbf{k}[X]) \subset \mathcal{I}\big(\mathbf{k}[Y])$.

\par
In \cite{ENSS}, the authors have shown that if $X$ is a finite quandle having a subquandle $Y$ with more than one element such that $|Y|$ is invertible in $\mathbf{k}$. Then $\mathbf{k}[X]$ has a non-trivial idempotent. Based on this, we have the following proposition whose proof is straightforward.

\begin{prop}
    Let $X$ be a finite quandle of odd order. Then $\sum_{x\in X} e_x$ is an idempotent in $\mathbb{Z}_2[X]$.
\end{prop}

Note that for some quandles $X$, the set ($\mathcal{I}\big(\mathbf{k}[X]), \cdot)$ with multiplication may not always be a quandle.  The following example shows that the set of idempotents over $\mathbb{Z}_2$ does not satisfy the condition that right multiplication by a fixed element being invertible and thus it is not a quandle.  First, we need the following notation. 

\noindent
{\bf Notation:}  Through the rest of the article, the notation $C[i,j]$ will mean 
the $j$-th connected quandle of order $i$ (for more details see \cite{rig}).

\begin{example}


Let $X=C[6,1]$ be the first connected quandle in the list of connected quandle \cite{rig} of order six.  As a set $ X=\{1, 2, 3,..., 6\}$ and its binary operation is given by the following cayley table:

\begin{center}
    \begin{tabular}{c|c c c c c c}
    
        $\ast$ & 1 & 2 & 3 & 4 & 5 & 6 \\
         \hline
        1 & 1 & 1 & 5 & 6 & 3 &4 \\
        
        2 & 2 & 2 & 6 & 5 & 4 & 3 \\
        
        3 & 5 & 6 & 3 & 3 & 1 & 2 \\
         
        4 & 6 & 5 & 4 & 4 & 2 & 1 \\
         
        5 & 3 & 4 & 1 & 2 & 5 & 5\\
        
        6 & 4 & 3 & 2 & 1 & 6 & 6
    \label{tab:c61}
    \end{tabular}
\end{center}
From the table of the quandle $X$, we see that
$e_1+e_3+e_5$ is an idempotent element in $Y=\mathcal{I}\big(\mathbb{Z}_2[X])$.  Now since $e_1(e_1+e_3+e_5)=e_1+e_3+e_5$ and $e_3(e_1+e_3+e_5)=e_1+e_3+e_5$, the right multiplication by the element $e_1+e_3+e_5$ fails to be injective and thus $Y$ is not a quandle. 
The article \cite{ENS} provides a table of quandles up to order 5 with their idempotents computed with coefficients in $\mathbb{Z}$ and $\mathbb{Z}_2$ and states the cases when such a set of idempotents is a quandle.
\end{example}

    

\section{Distinguishing knots of 12 crossings with cocycle invariants and idempotents in quandle rings}\label{Sec4} 


  It is known \cite{CJKLS} that the coloring invariant of link is weaker than the quandle 2-cocycle invariant of the link.
 Below is an example of classes of knots which were not distinguished by
 coloring \cite{CESY}. However, we were able to distinguish them by using the 2-cocycle of $C[12,3]$,  $\mathcal{I}\big( \mathbb{Z}_2[C[12,3]]\big) $ and $C[13,4]$.
\begin{example}
Let $X= C[12,3]$ be the third connected quandle \cite{rig} of order 12.  As a set $X=\{1,2,...,12\}$.  Its quandle operation is given in terms of right multiplications as follows: \newline

$\begin{array}{cc}
 S_1= (2~12~5~10~11)(3~8~6~7~4)   & S_2= (1~11~7~4~12)(3~5~10~6~9) \\ S_3=(1~2~7~6~10)(4~9~8~5~12) &
S_4=(1~11~6~8~5)(2~7~9~3~12) \\   S_5=(1~12~3~8~10)(2~4~9~6~11) & S_6=(1~5~3~4~2)(7~11~10~8~9)\\
S_7=(1~10~8~3~12)(2~11~6~9~4) & S_8=(1~12~5~7~11)(3~9~6~10~5) \\
 ~S_9=(2~11~10~5~12)(3~4~7~6~8) &
 S_{10}=(1~5~8~6~11)(2~12~3~9~7) \\ S_{11}=(1~10~6~7~2)(4~12~5~8~9) & ~S_{12}=(1~2~4~3~5)(7~9~8~10~11).
\end{array} $

\noindent
Using Maple software, we obtained the following $2$-cocycle with coefficients in $\mathbb{Z}_2$.  The map $\phi: X \times X \to \mathbb{Z}_2$ is given explicitly by \\
$\begin{array}{ccccc}
  \phi(3,2)=1, & \phi(3,4)=1, & \phi(4,7)=1, & \phi(4,11)=1, & \phi(6,5)=1, \\ \phi(6,8)=1, &
  \phi(7,6)=1, & \phi(8,3)= 1, & \phi(8,12)=1, & \phi(9,2)=1, \\ \phi(9,3)=1, & \phi(9,4)=1, &
  \phi(9,5)= 1, & \phi(9,6)=1, & \phi(9,7)=1, \\
  \phi(9,8)=1, & \phi(9,10)= 1, & \phi(9,11)=1, &
   \phi(9,12)=1 \\
  
 \end{array}$ \\
  and $\phi(x,y)=0$ for all other  $x, y \in Y$. \\ 
Note that this 2-cocycle is not a coboundary since the value of the quandle cocycle invariant of the knot $12n_{368}$ is given by  $ \Phi_{(C[12,3], \phi)}(12n_{368})= 40+32u $. \\
The 2-cocycle invariant $ \Phi_{(X, \phi)}(K)$ of the knots $K \in \{ 9_{13}, 9_{14}, 9_{16}, 9_{20}, 9_{23}, 9_{24}, 10_{123}, 12n_{0572}, \\12n_{0576},12n_{0578}, 12n_{0580} \}$  has value $72$. \\
To distinguish these knots further, we use the following quandle  $Y=\mathcal{I}\big(\mathbb{Z}_{2}  [X]\big)$.  As a set, we write  $Y=\{1,2,...,24\}$ and we give its quandle structure by listing its right multiplication given below: \newline
 $S_1=S_{13}= (2~12~5~10~11)(3~8~6~7~4)(14~24~17~22~23)(15~20~18~19~16)\\
S_2=S_{14}= (1~11~7~4~12)(3~5~10~6~9)(13~23~19~16~24)(15~17~22~18~21)\\
S_3=S_{15}=(1~2~7~6~10)(4~9~8~5~12)(13~14~19~18~22)(16~21~20~17~24)\\ 
S_4=S_{16}=(1~11~6~8~5)(2~7~9~3~12)(13~23~18~20~17)(14~19~21~15~16)$\\
$S_5=S_{17}=(1~12~3~8~10)(2~4~9~6~11)(13~24~15~20~22)(14~16~21~18~23)\\
S_6=S_{18}=(1~5~3~4~2)(7~11~10~8~9)(13~17~15~16~14)(19~23~22~20~18) \\
S_7=S_{19}=(1~10~8~3~12)(2~11~6~9~4)(13~22~20~15~24)(14~23~18~21~16) \\
S_8=S_{20}=(1~12~5~7~11)(3~9~6~10~5)(13~24~16~19~23)(15~21~18~22~17) \\
S_9=S_{21}=(2~11~10~5~12)(3~4~7~6~8)(14~23~22~17~24)(15~16~19~18~20) \\
S_{10}=S_{22}=(1~5~8~6~11)(2~12~3~9~7)(13~17~20~18~23)(14~24~15~21~19) \\
S_{11}=S_{23}=(1~10~6~7~2)(4~12~5~8~9)(13~22~18~19~14)(16~24~17~20~21) \\
S_{12}=S_{24}=(1~2~4~3~5)(7~9~8~10~11)(13~14~16~15~17)(19~21~20~22~23) \\
$
with 2-cocycle map $\psi: Y \times Y \to \mathbb{Z}_2$ given by \newline
$\begin{array}{ccccc}
~~~\psi(3,2)=1, & \psi(3,4)=1, & \psi(3,16)=1, & \psi(4,7)=1, &
\psi(4,11)=1, \\ 
\psi(4,19)=1, & \psi(4,23)=1, & \psi(7,6)=1, & \psi(7,10)=1, & \psi(7,18)=1, \\
 \psi(7,22)=1, & \psi(8,3)=1, & \psi(8,12)=1, & \psi(8,15)=1, & \psi(8,24)=1, \\
 \psi(9,2)=1, &\psi(9,3)=1, & \psi(9,4)=1, & ~\psi(9,5)=1, & \psi(9,6)=1, \\
 \psi(9,7)=1, & \psi(9,8)=1, & \psi(9,10)=1, & \psi(9,11)=1, &
~~\psi(9,12)=1, \\
\psi(9,14)=1, & \psi(9,15)=1, &\psi(9,20)=1, & \psi(9,22)=1, &\psi(9,23)=1, \\
 ~\psi(9,24)=1, & \text{and}~ \psi(x,y)=0 & \text{for all other}~x, y \in Y. \\
\end{array}$  \\ 
This further breaks down the set of knots into the following partition $\{  9_{13}, 9_{14}, 9_{16}, 9_{20}, 9_{23}, 9_{24} \} \sqcup \{10_{123}, 12n_{0572} \} \sqcup \{12n_{0576}\} \sqcup \{12n_{0578}\} \sqcup \{ 12n_{0580} \}$ since the cocycle invariants for each partition are respectively $144, 106+38u, 58+86u, 120+24u$ and $64+80u$.\\
To completely distinguish all the knots, we use the following quandle $C[13,4]$.  As a set we denote it by $W=\{ 1,2,3,...,13\}$.  Its quandle operation is given in terms of right multiplications by
\newline
$S_1= (2~9~13~6)(3~4~12~11)(5~7~10~8)~~
S_2= (1~7~3~10)(4~5~13~12)(6~8~11~9)\\
S_3=(1~13~5~6)(2~8~4~11)(7~9~12~10)~~
S_4=(1~6~7~2)(3~9~5~12)(8~10~13~11)\\
S_5=(1~12~9~11)(2~7~8~3)(4~10~6~13)~~
S_6=(1~5~11~7)(2~13~10~12)(3~8~9~4) \\
S_7=(1~11~13~3)(2~6~12~8)(4~9~10~5)~~
S_8=(1~4~2~12)(3~7~13~9)(5~10~11~6)\\
S_9=(1~10~4~8)(2~5~3~13)(6~11~12~7) ~~
S_{10}=(1~3~6~4)(2~11~5~9)(7~12~13~8)\\
S_{11}=(1~9~8~13)(2~4~7~5)(3~12~6~10)~~
S_{12}=(1~2~10~9)(3~5~8~6)(4~13~7~11)\\
S_{13}=(1~8~12~5)(2~3~11~10)(4~6~9~7)$. \\
Now we consider the following $2$-cocycle map $\vartheta: W \times W \to \mathbb{Z}_3$ given by \\
$\begin{array}{ccccc}
~~\vartheta(1,12)=1,  & \vartheta(2,10)=1, & \vartheta(3,8)=1, & \vartheta(4,6)=1, & \vartheta(5,4)=1, \\
~\vartheta(6,2)=1, &  \vartheta(7,13)=1, & \vartheta(8,11)=1, & \vartheta(9,1)=2, & \vartheta(9,2)=2, \\
~\vartheta(9,3)=2,& \vartheta(9,4)=2, &  \vartheta(9,5)=2, &  \vartheta(9,6)=2, & \vartheta(9,7)=2, \\
~\vartheta(9,8)=2, & \vartheta(9,10)=2, &  \vartheta(9,11)=2, &  \vartheta(9,12)=2, & \vartheta(9,13)=2, \\
~~\vartheta(10,7)=1, & \vartheta(11,5)=1, & \vartheta(12,3)=1, &\vartheta(13,1)=1, & \text{and} \\
~\vartheta(x,y)=0 & \text{for all other}~ x, y \in W.\\
\end{array}$  \\
 
Finally, we are able to distinguish all the above knots in the following table:

\begin{center}
    \begin{tabular}{|c|c|c|c|}
\hline
   $K$ & $\Phi_{(X, \phi)}(K)$ & $\Psi_{(Y, \psi)}(K) $ & $\Theta_{(W, \vartheta)}(K)$ \\
   \hline
    $9_{13}$ & 72 & 144 & $ 11u^2+42u+52$\\
    \cline{4-4}
    $9_{14}$ & 72 & 144 & $ 20u^2+5u+14 $\\
    \cline{4-4}
    $9_{16}$ & 72 & 144 & $ 69+14u $\\
    \cline{4-4}
    $9_{20}$ & 72 & 144 & $ 7u^2+13u+20 $\\
    \cline{4-4}
    $9_{23}$ & 72 & 144 & $ 13 $\\
    \cline{4-4}
    $9_{24}$ & 72 & 144 & $ 45u^2+37u+13 $\\
     \cline{3-4}
    $10_{123}$ & 72 & $106+38u$ & $3u^2+20u+17$ \\
    \cline{4-4}
    $12_{n0572}$ & 72 & $106+38u$ & $ 100u^2+70u+11 $ \\
    \cline{3-4}
    $12_{n0576}$ & 72 & $ 58+86u $ & $ 20u^2+135u+10 $\\
    \cline{3-4}
    $12_{n0578}$ & 72 & $ 120+24u $ & $ 16u^2+104u+13 $ \\
    \cline{3-4}
    $12_{n0580}$ & 72 & $ 64+80u $ & $ 78u^2+90u+11 $\\
   \hline
\end{tabular} \\
\end{center}
In a similar manner, we present a table that distinguish some non-alternating knots of 12 crossings.
\begin{center}
	\begin{tabular}{|c|c|c|c|c|}
		\hline
		$K$ & $ \mathcal{C}_X(K) $ & $ \Phi_{(X, \phi)}(K) $ & $ \Psi_{(Y, \psi)}(K) $ \\
	\hline
	$ 12n_{0573} $ & 132  & $ 68+64u $ & $ 136+128u $\\
	\cline{3-4}
	$ 12n_{0575} $ & 132  & $ 48+84u $ & $ 172+92u $\\
	\cline{4-4}
	$ 12n_{0577} $ & 132  & $ 48+84u $ & $ 144+120u $\\
	\cline{3-4}
	$ 12n_{0579} $ & 132  & $ 76+56u $ & $ 96+168u $\\
    \hline
	$ 12n_{0581} $ & 192  & $ 94+98u $ & $ 196+188u $ \\
    \cline{4-4}
        $ 12n_{0594} $ & 192  & $ 94+98u $ &  $ 240+ 144u $          \\
	\hline
	$ 12n_{0574} $ & 312  & $ 72+240u $ & $ 320+304u $ \\
  \cline{4-4}
        $ 12n_{0737} $ & 312  &  $ 72+240u $ & $ 324+ 300u $           \\
	\hline
	
	\end{tabular}
		\end{center}
In conclusion, the triplet ($\Phi_{(X, \phi)}(K), \Psi_{(Y, \phi)}(K), \Theta_{(W, \vartheta)}(K))$ distinguishes all the above mentioned knots completely.  Following this strategy we were able to distinguish all knots up to 12 crossings using $13$ quandles listed in the (\nameref{sec:appendix} below).  

\end{example}

 Based on the above computation, we have the following conjecture.
\begin{conjecture}
\label{conj1}
   Let X be a quandle and $\phi: X\times X \to A$ be a 2-cocycle.  Let $Y=\mathcal{I}(\mathbb{Z}_2[X])$ be the set of idempotents in the quandle ring $\mathbb{Z}_2[X]$ such that Y is a quandle and $\psi: Y \times Y \to A$  be a $2$-cocycle. Then the cocycle invariant $\Psi_{(Y,\psi)}(K)$ is an enhancement of the cocycle invariant $\Phi_{(X,\phi)}(K)$ for all prime oriented knots up to 12 crossings.

\end{conjecture}
This can further be extended to give the following generalized conjecture
\begin{conjecture}
\label{conj2}
    There exists a finite sequence of quandles $(X_1, X_2, X_3,...,X_k)$ such that \newline \small{$\Psi(K)=\big(\Phi_{(X_1, \phi_1)}(K), ...,\Phi_{(X_k, \phi_k)}(K),\Phi_{(\mathcal{I}(\mathbb{Z}_2[X_1]), \psi_1)}(K), ...,\Phi_{(\mathcal{I}(\mathbb{Z}_2[X_k]), \psi_k)}(K)\big)$ }is an invariant.  In other words, \\
    $\Phi(K)=\Phi(K')$ if and only if $K=K'$ for all $K$ in the list of knots up to 12 crossings.
\end{conjecture}




\begin{example}
For a knot $K$ let $m(K)$ denote the mirror image of $K$. We say that  a knot $K$ is {\it positive amphicheiral} if  $K = m(K)$. In this example, we give two knots which are not distinguished from their mirror image by the \emph{Jones} polynomial or by the quandle coloring, but we are able to distinguish them using the pair of 2-cocycle invariants $(\Phi_{(X, \phi)}(K), \Phi_{(Y, \psi)}(K))$ where $X = C[12,3]$ and $ Y=\mathcal{I}\big(\mathbb{Z}_{2}[C[12,3]]\big)$.\\
\begin{table}[ht]
\begin{tabular}{|c|c|c|c|c|c|}
\hline
   $K$  & $\mathcal{C}_{X}(K)$ & Jones Polynomial &  $\big(\Phi_{(X,\phi)}(K), \Psi_{(Y, \psi)}(K)\big)$ \\
\hline
   $9_{42}$  & $ 24 $ & $t^{-3}-t^{-2}+ t^{-1}-1+$ & $ (24, 48) $\\
   \cline{1-1}
   \cline{4-4}
   $m(9_{42})$ & $ 24 $  & $t-t^2+ t^3$ & $ (24, 32u+16 ) $\\
   \cline{1-1}
   \cline{3-4}
   $ 12a_{669}$& $24$ & $ -t^{-6}+ 2t^{-5}-4t^{-4}+  6t^{-3}-7t^{-2}+ 9t^{-1}-9+$ &$ (24, 48 ) $ \\
   \cline{1-1}
   \cline{4-4}
   $ m(12a_{669})$& $24$ & $ 9t-7t^2+ 6t^3-4t^4+ 2t^5-t^6$ &$ (24, 12u+36 ) $\\
   \hline
\end{tabular}
    
    \label{tab:my_label}
\end{table}

\end{example}



 
\section{Distinguishing knots of 13 crossings with cocycle invariants and idempotents in quandle rings}\label{Sec6}

We extended our computation towards 13 crossing knots for both alternating and non-alternating to support our Conjecture \ref{conj2}. It turns out that we needed 24 quandles to distinguish these knots up to 13 crossings (see the \nameref{sec:appendix} below for the list of the quandles). For example, consider the quandle $X=C[12,3]$ and $Y= \mathcal{I}\big(\mathbb{Z}_{2}[X]\big)$ (see \nameref{sec:appendix} below). We compute the 2-cocycle invariants $\Phi_{(X,\phi)}(K)$ and $\Psi_{(Y,\psi )}(K)$ using these two quandles. Now, we iterate the process of taking idempotents and consider the set $W= \mathcal{I}\big(\mathbb{Z}_{2}[\mathcal{I}\big(\mathbb{Z}_{2}[X]\big)]\big)$. We find that this set $W$ also forms a quandle. Additionally, the 2-cocycle invariant $\Theta_{(W, \vartheta)}(K)$ from $W$ is stronger than the previous 2-cocycle invariant $\Phi_{(X,\phi)}(K)$ and $\Psi_{(Y, \psi)}(K)$ from $X$ and $Y$ respectively. 
The following table further supports the claim:
\begin{center}

    \begin{tabular}{|c|c|c|c|c|}
\hline
   $K$  &  $\Phi_{(X,\phi)}(K)$ & $\Psi_{(Y,\psi )}(K)$  & $\Theta_{(W,\theta)}(K)$ \\
\hline
   $ 13_{120} $  &  $ 132 $ & $ 264 $ & $ 528 $\\
  \cline{4-4}
   $ 13_{482} $  &  $ 132 $ & $ 264 $ & $ 280+248u $\\
   \hline
   $ 13_{484} $  &  $ 48+84u $ & $ 144+120u $ & $ 252+276u $\\
   \cline{4-4}
   $ 13_{485} $  &  $ 48+84u $ & $ 144+120u $ & $ 240+288u $\\
  \cline{4-4}
    $ 13_{1596} $  &  $ 48+84u $ & $ 144+120u $ & $ 360+168u $\\
    \hline
\end{tabular}
\end{center}
The article \cite{CESY} defines some similarity of quandles. Precisely, for any two quandles $X_1 , X_2 $ and a family of knots $\mathcal{K}$, we say $X_1 \approx X_2$ if $ \mathcal{C}_{X_1}(K)=\mathcal{C}_{X_2}(K)$ for every $K \in \mathcal{K}$. 
\par In a similar manner we introduce a more general similarity of quandles using 2-cocycle invariant of knots; i.e.
for any two quandles and their respective 2-cocycles given by $(X_1,\phi), (X_2, \psi)$ and a family of knots $\mathcal{K}$, we say $X_1 \sim X_2$ if $ \Phi_{(X_1, \phi)}(K)=\Psi_{(X_2, \psi)}(K)$ for every $K \in \mathcal{K}$. 

Based on this, we have the following observation for all 12965 prime oriented knots up to 13 crossings and the 24 quandles in our computation:

\begin{itemize}
    \item There are a total of 3520 classes of $\sim$ consisting more than one quandle for knots up to 13 crossings. For example:
    \begin{itemize}
    \item  $C[12,3] \sim C[12,6]$ for $K\in \mathcal{K}=\{ 9_2, 9_3, 9_4, 9_7, 9_9, 9_{12}, 9_{13}, 11a_{172},11a_{190}, 11a_{191},\\ 12n_{0370}, 12n_{0371},  12n_{0373},  12n_{376}, 13_{3108} \}$.
        \item $ C[13,7] \sim C[13,10]$ for $ K\in \mathcal{K}= \{7_4, 7_6, 7_7, 8_7, 8_9, 8_{13}, 8_{17}, 8_{18}, 10_{46}, 10_{47}, 10_{48},\\ 10_{49}, 10_{50}, 10_{51}, 10_{52}, 11a_{151}, 11a_{152}, 11a_{171}, 13_{528}, 13_{3109}, 13_{9089}, 13_{9090}\}$.
    \end{itemize}
      \item From the $3520$ classes, there are $1460$ classes containing more than two quandles for knots up to 13 crossings. 
      \begin{itemize}
           \item $C[12,3] \sim C[12,4]\sim C[12,6]$ for $K\in \mathcal{K}=\{ 9_3, 9_5, 9_7, 9_8, 9_9, 12n_{0370}, 12n_{0371},\\ 12n_{0373}, 12n_{376}\}$  . 
          \item $C[16,3]\sim C[16,4] \sim \mathcal{I}\big(\mathbb{Z}_2[C[8,1]]\big)$ for $K\in \mathcal{K}=\{ 12n_{370}, 12n_{371}, 12n_{372}, 12n_{373},\\ 12n_{374}, 12n_{376},  13_{4039}, 13_{535}, 13_{544} \}$ 
      \end{itemize}
        \end{itemize}

\section{A Brief Description of the Algorithm}\label{Sec7}

In this section, we give a brief description of the algorithm. The algorithm was inspired by 
\url{http://shell.cas.usf.edu/~saito/Maple/}. The algorithm has three steps  
\begin{itemize}
    \item {\bf Step 1.} Given a quandle X, we check when the set of idempotents $\mathcal{I}\big(\mathbb{Z}_2[X]\big)$ is a quandle.
    \item {\bf Step 2.} Using a quandle $X$ and its idempotent $\mathcal{I}\big(\mathbb{Z}_2[X]\big)$ over $\mathbb{Z}_2$, we calculate the colorings of a given knot $K$ using its braid representation.
    \item {\bf Step 3.} After obtaining the coloring, we calculate the State Sum Invariant of the knot $K$ using both $X$ and $\mathcal{I}\big(\mathbb{Z}_2[X]\big)$ .
\end{itemize}
Below we give a quick description of how to find the coloring of a knot $K$ and the State Sum Invariant of $K$.
\begin{itemize}
\item {\bf Finding the colorings of a knot $K$:} 
Since every knot is the closure of a braid, to obtain a coloring of the knot, we first color the braid as follows.
\begin{enumerate}
    \item Let $m$ be the braid index of the knot $K$.
    \item 
    Let $\vec{x}=(x_1, x_2, x_3,...,x_i,x_{i+1},...,x_m) \in X^m$ be the top color of the braid.
    \item 
    At any $i^{th}$ crossing, if it is positive, then $\vec{x}$ becomes $(x_1, ...,x_{i+1},x_i\ast x_{i+1},...,x_m)$. if the $i^{th}$ crossing is negative, then $\vec{x}$ becomes $(x_1,...,x_{i+1}\Bar{\ast} x_i,x_i,...,x_m)$. Let $(y_1, y_2, ..., y_m)$ be the bottom vector of the braid.
    \item A solution of the system of equations $y_1=x_1,...,~y_m=x_m$ is a coloring of knot $K$ by the quandle $X$.  We abuse the notation and use $\vec{x}$ to denote this coloring.
\end{enumerate}

\item {\bf Calculating the State Sum Invariant of $K$}
\begin{enumerate}
    \item Let $\vec{x}=(x_1, x_2, x_3,...,x_i,x_{i+1},...,x_m) \in X^m$ be the top color of the braid, thus giving a color of $K$
    \item At a positive crossing with input colors $x_i$ and $ x_{i+1}$ , we assign the Boltzmann weight $\phi(x_i, x_{i+1})$ (refer to the left picture of Figure \ref{boltz}). Similarly at a negative crossing with output colour $x_i$ and $ x_{i+1}$ we assign the Boltzmann weight $-\phi(x_i, x_{i+1})$ (refer to the right picture of Figure \ref{boltz}). 
    \item For each fixed coloring, we compute the product $\prod_\tau \phi(x,y)^{\epsilon}$ over all crossings. The result is an element of group $A$. To obtain the State Sum Invariant of $K$, we sum over all possible coloring obtaining the element $\Phi(K)= \sum_{\mathcal{C}} \prod_\tau \phi(x,y)^{\epsilon}$ of the integral group ring $A$.
\end{enumerate}
\end{itemize}
\section*{Acknowledgement} 
ME was partially supported by Simons Foundation collaboration grant 712462.
The authors would like to thank Masahico Saito for fruitful discussions which improved the paper. DS thanks Manpreet Singh for his help with the Python Programming.
\\

\vspace{3ex}
\appendix
\section*{{Appendix}}\label{sec:appendix}

In this Appendix we provide the list of quandles and some of their idempotent quandles 
used in this paper explicitly (a full list of the 24 quandles can be found at \url{https://swaindipali.com/wp-content/uploads/2024/02/list-of-quandles.txt} ).\\ 
{\bf Note:} The notation $C[i,j]$ stands for the $j$-th connected quandle of order $i$ (see \cite{rig}).  The right multiplications $S_k$ in the quandle are given by $S_k(l)=l*k$.  As permutations, right multiplications are written below as product of cycles.

\begin{itemize}

\small{
\item $C[8,1]$ \\
$S_1=S_2= (3~6~7)(4~5~8) ~~~~~
S_3=S_4=(1~8~6)(2~5~7) ~~~~~
S_5=S_6=(1~4~7)(2~3~8)\\
S_7=S_8=(1~5~3)(2~6~4)$ \\

\item $\mathcal{I}\big(\mathbb{Z}_2[C[8,1]] \big)$.  As a set $\mathcal{I}\big(\mathbb{Z}_2[C[8,1]] \big)=\{1, \ldots, 16\}$.\\
$S_1=S_2=S_{9}=S_{10}= (3~6~7)(4~5~8)(15~11~14)(12~13~16) \\
S_3=S_4=S_{11}=S_{12}=(1~8~6)(2~5~7)(13~10~15)(9~16~14) \\
S_5=S_6=S_{13}=S_{14}=(1~4~7)(2~3~8)(9~12~15)(10~11~16) \\
S_7=S_8=S_{15}=S_{16}=(1~5~3)(2~6~4)(9~13~11)(10~14~12)$ \\
\item $C[12,3]$ \\
$\begin{array}{cc}
  S_1= (2~12~5~10~11)(3~8~6~7~4)   & ~S_2= (1~11~7~4~12)(3~5~10~6~9) \\ 
  S_3=(1~2~7~6~10)(4~9~8~5~12) &   S_4=(1~11~6~8~5)(2~7~9~3~12)\\
  S_5=(1~12~3~8~10)(2~4~9~6~11) & S_6=(1~5~3~4~2)(7~11~10~8~9) \\
  ~S_7=(1~10~8~3~12)(2~11~6~9~4) & S_8=(1~12~5~7~11)(3~9~6~10~5) \\
    S_9=(2~11~10~5~12)(3~4~7~6~8) &
  S_{10}=(1~5~8~6~11)(2~12~3~9~7) \\
  S_{11}=(1~10~6~7~2)(4~12~5~8~9) & S_{12}=(1~2~4~3~5)(7~9~8~10~11).\\
\end{array} \\$

   \item $ \mathcal{I}\big(\mathbb{Z}_2[C[12,3]] \big) $.   As a set $\mathcal{I}\big(\mathbb{Z}_2[C[12,3]] \big)=\{1, \ldots, 24\}$. \\
    $S_1=S_{13}= (2~12~5~10~11)(3~8~6~7~4)(14~24~17~22~23)(15~20~18~19~16)\\
S_2=S_{14}= (1~11~7~4~12)(3~5~10~6~9)(13~23~19~16~24)(15~17~22~18~21)\\
S_3=S_{15}=(1~2~7~6~10)(4~9~8~5~12)(13~14~19~18~22)(16~21~20~17~24)\\ 
S_4=S_{16}=(1~11~6~8~5)(2~7~9~3~12)(13~23~18~20~17)(14~19~21~15~16)$\\
$S_5=S_{17}=(1~12~3~8~10)(2~4~9~6~11)(13~24~15~20~22)(14~16~21~18~23)\\
S_6=S_{18}=(1~5~3~4~2)(7~11~10~8~9)(13~17~15~16~14)(19~23~22~20~18) \\
S_7=S_{19}=(1~10~8~3~12)(2~11~6~9~4)(13~22~20~15~24)(14~23~18~21~16) \\
S_8=S_{20}=(1~12~5~7~11)(3~9~6~10~5)(13~24~16~19~23)(15~21~18~22~17) \\
S_9=S_{21}=(2~11~10~5~12)(3~4~7~6~8)(14~23~22~17~24)(15~16~19~18~20) \\
S_{10}=S_{22}=(1~5~8~6~11)(2~12~3~9~7)(13~17~20~18~23)(14~24~15~21~19) \\
S_{11}=S_{23}=(1~10~6~7~2)(4~12~5~8~9)(13~22~18~19~14)(16~24~17~20~21) \\
S_{12}=S_{24}=(1~2~4~3~5)(7~9~8~10~11)(13~14~16~15~17)(19~21~20~22~23) \\
$

\item $C[12,4]$\\
$S_1= (5~9)(2~3~4)(6~11~8~10~7~12) ~~~S_2=(6~10)(1~4~3)(5~12~7~9~8~11)\\
S_3=(7~11)(1~2~4)(5~10~8~9~6~12)$~~~
$S_4=(8~12)(1~3~2)(5~11~6~9~7~10)~~~~\\
S_5=(1~9)(6~7~8)(2~11~4~10~3~12)~~~S_6=(2~10)(5~8~7)(1~12~3~9~4~11)$ \\
$S_7=(3~11)(5~6~8)(1~10~4~9~2~12)~~~S_8=(4~12)(5~7~6)(1~11~2~9~3~10)\\
S_9=(1~5)(10~11~12)(2~7~4~6~3~8)$ 
$S_{10}=(2~6)(9~12~11)(1~8~3~5~4~\\
S_{11}=(3~7)(9~10~12)(1~6~4~5~2~8)~~~S_{12}=(4~8)(9~11~10)(1~7~2~5~3~6)$.\\
\item $C[12,6]$\\
$S_1=(3~4)(5~10)(6~9)(8~12)(7~11)  ~~~S_2=(3~4)(5~9)(6~10)(8~11)(7~12) \\
S_3=(1~2)(5~11)(6~12)(7~9)(8~10)$~~~ $S_4=(1~2)(5~12)(6~11)(7~10)(8~11)\\
S_5=(1~10)(2~9)(3~11)(4~12)(7~8)~~~S_6=(1~9)(2~10)(3~12)(4~11)(7~8)$ \\
$S_7=(1~11)(2~12)(3~9)(4~10)(5~6)~~~S_8=(1~12)(2~11)(3~10)(4~9)(5~6)\\
S_9=(1~6)(2~5)(3~7)(4~8)(11~12)$ 
$S_{10}=(1~5)(2~6)(3~8)(4~7)(11~12)\\
S_{11}=(1~7)(2~8)(3~5)(4~6)(9~10)~~~S_{12}=(1~8)(2~7)(3~6)(4~5)(9~10)$.\\

\item $C[13,4] $ \\
$S_1= (2~9~13~6)(3~4~12~11)(5~7~10~8) ~~
S_2= (1~7~3~10)(4~5~13~12)(6~8~11~9)\\
S_3=(1~13~5~6)(2~8~4~11)(7~9~12~10)~~~
S_4=(1~6~7~2)(3~9~5~12)(8~10~13~11)~~\\
S_5=(1~12~9~11)(2~7~8~3)(4~10~6~13)~~~
S_6=(1~5~11~7)(2~13~10~12)(3~8~9~4)\\
S_7=(1~11~13~3)(2~6~12~8)(4~9~10~5)~~
S_8=(1~4~2~12)(3~7~13~9)(5~10~11~6)\\
S_9=(1~10~4~8)(2~5~3~13)(6~11~12~7)~~~
S_{10}=(1~3~6~4)(2~11~5~9)(7~12~13~8)~~\\
S_{11}=(1~9~8~13)(2~4~7~5)(3~12~6~10)~~~
S_{12}=(1~2~10~9)(3~5~8~6)(4~13~7~11)~~~\\
S_{13}=(1~8~12~5)(2~3~11~10)(4~6~9~7)$. \\

\item $C[13,7]$\\
$S_1= (2~6~13~9)(3~11~12~4)(5~8~10~7) ~~
S_2= (1~10~3~7)(4~12~13~5)(6~9~8~11)~~~\\
S_3=(1~6~5~13)(2~11~4~8)(7~10~12~9)~~~
S_4=(1~2~7~6)(3~12~5~9)(8~11~13~10)~~\\
S_5=(1~11~9~12)(2~3~8~7)(4~13~6~10)~~~
S_6=(1~7~11~5)(2~!2~10~13)(3~4~9~8)\\
S_7=(1~3~13~11)(2~8~12~6)(4~5~10~9)~~
S_8=(1~12~2~4)(3~9~13~7)(5~6~11~10)~~~\\
S_9=(1~8~4~10)(2~13~3~5)(6~7~12~11)~~~~
S_{10}=(1~4~6~3)(2~9~5~11)(7~8~13~12)~~\\
S_{11}=(1~13~8~9)(2~5~7~4)(3~10~6~12)~~~
S_{12}=(1~9~10~2)(3~6~8~5)(4~11~7~13)~~~\\
S_{13}=(1~5~12~8)(2~10~11~3)(4~7~9~6)$. \\

\item $C[13,10]$\\
$S_1= (2~7~11~9~10~3~13~8~4~6~5~12) ~~
S_2= (1~9~5~7~6~13~3~8~12~10~11~4)~~~\\
S_3=(1~4~9~13~11~12~5~2~10~6~8~7)~~~
S_4=(1~12~13~6~3~11~7~9~8~2~5~10)~~\\
S_5=(1~7~4~12~8~10~9~3~6~11~2~13)~~~
S_6=(1~2~8~5~13~9~11~10~4~7~12~3)\\
S_7=(1~10~12~11~5~8~13~4~2~3~9~6)~~
S_8=(1~5~3~4~10~7~2~11~13~12~6~9)~~~\\
S_9=(1~13~7~10~2~6~4~5~11~8~3~12)~~~~
S_{10}=(1~8~11~3~7~5~6~12~9~4~13~2)~~\\
S_{11}=(1~3~2~9~12~4~8~6~7~13~10~5)~~
S_{12}=(1~11~6~2~4~3~10~13~5~9~7~8)~~~\\
S_{13}=(1~6~10~8~9~2~12~7~3~5~4~11)$. 

\item $C[16,3] $ \\
 $S_1= (2~3~5~9~16)(4~713~8~15)(6~11~12~10~14) ~~~~~
S_2= (1~4~6~10~15)(3~8~14~7~16)(5~12~11~9~13) ~~~\\
S_3=(1~7~11~14~4)(2~5~15~6~13)(8~9~10~12~16)~~~
S_4=(1~6~16~5~14)(2~8~12~13~3)(7~10~9~11~15)~~\\
S_5=(1~13~12~6~7)(2~15~16~14~10)(3~9~4~11~8)~~~
S_6=(1~16~15~13~9)(2~14~11~5~8)(3~12~7~4~10)\\
S_7=(1~11~2~9~6)(3~15~10~8~5)(4~13~14~16~12)~~~~
S_8=(1~10~5~2~12)(3~14~13~15~11)(4~16~9~7~6)~~~\\
S_9=(1~8~10~11~13)(2~6~13~3~4)(5~16~7~12~15)~~~~
S_{10}=(1~5~13~4~3)(2~7~9~12~14)(6~15~8~11~16)~~\\
S_{11}=(1~2~4~8~16)(3~6~12~9~15)(5~10~13~7~14)~~~
S_{12}=(1~3~7~15~2)(4~5~11~10~16)(6~9~14~8~13)~~~\\
S_{13}=(1~12~3~16~11)(2~10~7~8~6)(4~14~15~9~5)~~~
S_{14}=(1~9~8~7~5)(2~11~4~15~12)(3~13~16~10~6) \\
S_{15}=(1~14~9~3~10)(2~16~13~11~7)(4~12~5~6~8)~~~
S_{16}= (1~15~14~12~8)(2~13~10~4~9)(3~11~6~5~7)\\$
\item $C[16,4]$ \\
 $S_1=(2~13~4~9~3~5)(6~14~16~12~11~7)(8~10~15)  ~~~~~
S_2= (1~14~3~10~4~6)(5~13~15~11~12~8)(7~9~16) ~~~\\
S_3=(1~7~4~15~2~11)(5~8~16~14~10~9)(6~12~13)~~~
S_4=(1~12~2~8~3~16)(6~7~15~13~9~10)(5~11~14)~~\\
S_5=(1~6~9~8~13~7)(2~10~12~16~15~3)(4~14~11)~~~
S_6=(1~9~11~15~16~4)(2~5~10~7~14~8)(3~13~12)\\
S_7=(1~4~12~10~14~13)(3~8~11~6~15~5)(2~16~9)~~
S_8=(2~3~11~9~13~14)(4~7~12~5~16~6)(1~15~10)~~~\\
S_9=(1~11~13~10~5~12)(3~15~14~6~8~4)(2~7~16)~~~~
S_{10}=(2~12~14~9~6~11)(3~4~16~13~5~7)(1~8~15)~~\\
S_{11}=(1~13~16~8~6~2)(3~9~15~12~7~10)(4~5~14)~~~~
S_{12}=(1~2~14~15~7~5)(4~10~16~11~8~9)(3~16~13)~~~\\
S_{13}=(1~16~5~15~9~14)(2~4~8~7~11~10)(3~12~6)~~~
S_{14}=(1~3~7~8~12~9)(2~15~6~16~10~13)(4~11~5) \\
S_{15}=(2~6~5~9~12~4)(3~14~7~13~11~16)(1~10~8)~~
S_{16}=(1~5~6~10~11~3)(4~13~8~14~12~15)(2~9~7) \\$\\

\item $ C[16,8] $ \\
$S_1= (2~3~5~9~10~12~16~8~15~6~11~14~4~7~13) ~~~
S_2= (1~4~6~10~9~11~15~7~16~5~12~13~3~8~14~1) ~~~\\
S_3=(1~7~11~12~10~14~6~13~8~9~16~2~5~15~4)~~~
S_4=(1~6~16~3~2~8~12~11~9~13~5~14~7~10~15)~~\\
S_5=(1~13~14~16~12~4~11~2~15~10~8~3~9~6~7)~~~
S_6=(1~16~9~7~4~10~5~8~2~14~13~15~11~3~12)\\
S_7=(1~11~8~5~3~15~16~14~10~2~9~4~13~12~6)~~~~
S_8=(1~10~3~14~11~5~2~12~7~6~4~16~15~13~9)~~~\\
S_9=(1~2~4~8~16~7~14~3~6~12~15~5~10~11~13)~~~~
S_{10}=(1~3~7~15~8~13~4~5~11~16~6~9~12~14)~~\\
S_{11}=(1~8~10~13~7~12~9~15~3~4~2~6~14~5~16)~~~
S_{12}=(1~5~13~6~15~2~7~9~14~8~11~10~16~4~3)~~~\\
S_{13}=(1~14~15~9~5~6~8~4~12~3~10~7~2~16~11)~~~
S_{14}=(1~15~12~2~13~16~10~6~5~7~3~11~4~9~8) \\
S_{15}=(1~12~5~4~14~9~3~16~13~11~7~8~6~2~10)~~~
S_{16}=(1~9~2~11~6~3~13~10~4~15~14~12~8~7~5)$ \\
}
\end{itemize}

\end{document}